%% file: pressb.tex
\documentclass[10pt]{article}
\usepackage{amssymb}

\def\R{\mathbb{R}}
\def\C{\mathbb{C}}
\def\T{\mathbb{T}}
\def\D{\mathbb{D}}

\newtheorem{prop}{Proposition}

\newtheorem{theoreme}{Theorem}
\newtheorem{coro}{Corollary}
\newenvironment{dem}{\bigskip\noindent{\bf Proof} :}%

\newtheorem{exmp}{Example} 

\newenvironment{rem}{\medskip\noindent {\bf Remark \small} : }%

{\par}
{\par}

\begin{document}
\sf{
\begin{center}
\huge{\bf Asymptotic behavior of Patil's approximants in Hardy
spaces :\\ The real case }
\end{center}
\begin{center}
\vspace{1cm}
\large{\Large By Gomari Buanani.Naufal}\\
\end{center}
\begin{center}
\large{\Large Universit\'e Claude Bernard-Lyon Institut Girard
desargues. 43, Bd du 11 novembre 1918, 69622 Villeurbanne, France.
  gomari@desargues.univ-lyon1.fr}
\end{center}
\vspace{25 mm}
\large{\Large\bf Abstract} :\\
 
\textit{In this paper we consider a robust identification problem for
a linear dynamical control system with limited-frequency intervals. In
mathematical terms, this is the problem of recovering functions in
Hardy spaces.
Our purpose is to bound Patil's approximants in the upper half plane case, out of a
bounded real interval $I$. To this end, we deal with 
 residu techniques and  give a class of functions to provide
boundedness of these approximants on the complement of this
interval.}\\

\vspace{15 mm} 

\vspace{15 mm}
 

\pagebreak

\section{Introduction}

Let $ \D $ be the open unit disc and $ \T $ its boundary. $ H^{p}(\D) $
denotes the Hardy space of analytic functions on $ \D $.\\
In 1972, Patil in [9] has given an algorithm to recapture a $
H^{p}(\D) $, $1<p<\infty$
function F from its values on E, a positive Lebesgue
measure subset of $ \T $. Let $ f $ be the boundary function  of $ F $ and $g$ its restriction
 on $E$.\\
In fact, using Toeplitz operator techniques he has found a sequence of functions $ g_{\lambda} $ uniformly converging
 on compact subsets of $ \D $ to $ F $ and also strongly in $ H^{p}(\D)
$.\\
 The same approximants had been  obtained thanks to a Carlman's fruitful idea by constructing a ``quenching'' function,
enabling us to eliminate in the Cauchy formula, integration over
$\T\setminus E$. In fact, we first construct an auxilary function
$\varphi\in H^\infty (\D)$ satisfying :
\begin{enumerate}
\item $|\varphi (\xi)|=1$ ae on $\T\setminus E$.
\item $|\varphi (\xi)| > 1$ in $\D$.
\end{enumerate}
To do so, we solve a suitable Dirichlet problem. If $u$ is the bounded
harmonic solution in $\D$  such that :
$$
u(x,y)= \left\{ \begin{array}{ll}
           1\ \ \mbox{a.e on}\ E\\
           0 \ \ \mbox{on}\ \T\setminus E.
\end{array} \right.
$$
then we set $\varphi (z)=e^{u+iv}$ where $v$ is the harmonic function
conjugate with $u$. Thus, we obtain the following formula :
$$
f(z)=\lim_{m\rightarrow+\infty}\frac{1}{2\pi i}\int_{E}\left
(\frac{\varphi(\xi)}{\varphi(z)}\right)^m\frac{f(\xi)}{\xi-z}d\xi
$$

L.Baratchart, J.Leblond and J.R Partington have exhibited the same
approximants applying optimization methods [1,2]. Their work was initiated
by them and D.Alpay in the $H^2$ case, the parametre here is  a Lagrange multiplier in [1,2].

Mukherjee in [11], has dealt with the upper half plane case by exhibiting
the corresponding sequences $ h_{\lambda} $ via the Cayley transform
and $ g_{\lambda} $ using Wiener-Hopf operator. He has shown that results obtained
by Patil remain valid in this case.\\
We  wonder if in the case of the upper half plane we can hope better approximation of $ f $ out of $ I$ on
the boundary, where $I$ is a real bounded interval.\\
An open question is the almost everywhere convergence of $g_\lambda$ to $f$ on
$\T\setminus I$. 

Our aim is to show that in the upper half plane, under appropriate conditions,
$g_\lambda(x)=O(1)$ with respect to $\lambda$ out of $I$. We show also
that even if $g$ is the trace of a $H^2$ function on $I$, the
conjecture is false.   

We devide mainly our work into four sections. The first one is the
introduction. The second one
deals with nontangential limits of $ g_{\lambda} $ and $ h_{\lambda} $
 according to whether the interval is symmetrical around 0 or not. In the third
part, we give asymptotic behavior of $ g_{\lambda} $ with respect to $\lambda$ and exhibit a class
of functions such that the trace of $g_\lambda$ remains bounded 
 on $\R\setminus I$. Examples are given
in the fourth section.\\

\large{\bf Notations}:

\begin{enumerate}
\item $\C_+$ the open upper half plane.
\item $\chi _{I}$  the characteristic function on $I$.
\item p.v $\int$  the principal value of the integral.
\item $O(f)$  a function not exceeding f with respect to $\lambda$
when $\lambda \rightarrow +\infty$.
\end{enumerate}

\section{Boundary value functions}

Let recall the two main theorems of recovering functions in the both
cases $ \D $  and $ \C_{+} $. \\

\begin{theoreme}([9, thm 1])\\
Let $ E \subset \T$ with m($E$) $ >$ 0. Suppose that $g$ is the restriction of f
to $E$. For each $\lambda > 0$ define analytic functions $ h_{\lambda},\
g_{\lambda}$ on $\D$ by :

$$g_{\lambda}(z)=\lambda h_{\lambda}(z)\frac{1}{2\pi i}\int_K \frac
{\overline{h}_{\lambda}(w)g(w)}{w-z}dw \ z \in \D,$$

$$h_{\lambda}(z)=\exp \left (-\frac{1}{4\pi}\ln
(1+\lambda)\int_{K}\frac{e^{ i\theta}+z}{e^{ i\theta} -z}d\theta\right
) \    z \in \D.$$ 

$\mbox{Then as }\ \lambda\rightarrow\infty\ ;\ g_{\lambda}\rightarrow f\
\mbox{unifomly on compact subsets of $\D$. Moreover for}$ 
$1 < p <\infty\ \mbox{we also have}\ ||g_{\lambda}- f||_{p}\rightarrow 0\ as\
\lambda\rightarrow \infty.$
\end{theoreme}

\begin{theoreme}([11, thm 1])\\
Let $ I \subset \R $ with $ m(I)>0 $. Suppose that $F \in H^{2}(\C_+)$
and $g$ is the restriction of $f$ to $I$. For each $ \lambda > 0 $ we define
analytic functions on the upper half plane by :

$$h_{\lambda}(z)=\exp \left
(-\frac{1}{2\pi i}\ln(1+\lambda)\int_I\frac{1+tz}{(t-z)(1+t^{2})}dt)
\right )\Im(z)>0,$$
 
$$g_{\lambda}(z)=\lambda h_{\lambda}(z)\frac{1}{2\pi i}\int_I \frac
{\overline{h}_{\lambda}(t)g(t)}{t-z}dt \ \ \ \Im(z)>0.$$
then as $\lambda \rightarrow \infty $, $ g_{\lambda}\rightarrow F$
uniformly on compact subsets of the upper half plane. Moreover
$||g_{\lambda}- f||_{2}\rightarrow 0 $ as $\lambda \rightarrow \infty$.

\end{theoreme}

We treat two cases $ I=]-a,a[ $ where $ a > 0 $ and $ I=]a,b[ $ where
$0 \leq a < b$ .\\

\subsection{Symmetrical case}

Let  $ I=]-a,a[ $ $ a > 0 $, $g$ as in theorem 2.

\begin{prop}

If $ z=x+iy \in \C_{+}$ then $h_{\lambda}(z)\rightarrow h_{\lambda}(x)
$\ as \ $ y \rightarrow 0 $. Where $h_{\lambda}(x)=(1+\lambda
\chi_{I}(x))^{\frac{-1}{2}} e^{iG_{1}(x)}$, $G_{1}(x)=\displaystyle{\frac
{1}{2\pi}\ln(1+\lambda )\ln \left | \frac {a-x}{a+x}\right |}$.
\end{prop}

\begin{rem}
One can find this proposition in [8]. We give a correct proof of this.
\end{rem}

\begin{dem}

We know that 
$$h_{\lambda}(z)=\exp
\left (-\frac{1}{2\pi i}\ln(1+\lambda)\int_{-a}^{a}\frac{1+tz}{(t-z)(1+t^{2})}dt\right
)
for\ \Im(z) > 0.$$

We wish $\displaystyle\lim_{y\rightarrow 0} h_{\lambda}(z)$. We have 
$$
\int_{-a}^{a}\frac{1+tz}{(t-z)(1+t^{2})}dt=\int_{-a}^{a}\frac{1}{t-z}dt-\frac
{1}{2}\int_{-a}^{a}\frac{2t}{1+t^{2}}dt.
$$
On one hand $\displaystyle{\int_{-a}^{a}\frac{2t}{1+t^{2}}dt=0}$.\\

On the other hand, 
$$\int_{I}\frac{1}{t-z}dt=\int_{I}\frac {t-x}{(t-x)^{2}+y^{2}}dt +
i\int_{I}\frac {y}{(t-x)^{2}+y^{2}}dt.
 $$
It is easy to see that, 
$$\int_{I}\frac {t-x}{(t-x)^{2}+y^{2}}dt \rightarrow \ln\left | \frac
{a-x}{a+x}\right |\ \mbox {as}\ y \rightarrow 0 , $$
and
$$\ \int_{I}\frac {y}{(t-x)^{2}+y^{2}}dt \rightarrow \chi_{I}(x)\pi\ \mbox {as}\ y \rightarrow 0,  $$ 
so
$$
h_{\lambda}(z)\rightarrow \exp \left
(\frac {-1}{2\pi i}\ln(1+\lambda)\left [\ln\left | \frac {a-x}{a+x}\right
|+i\chi_{I}(x)\pi \right ] \right )\ \mbox {as}\ y\rightarrow 0.
 $$ 
Finally,
$$
h_{\lambda}(x)=(1+\lambda
\chi_{I}(x))^{\frac{-1}{2}} e^{iG_{1}(x)}, G_{1}(x)=\frac
{1}{2\pi}\ln(1+\lambda )\ln\left | \frac {a-x}{a+x}\right |.
 $$
\end{dem}

\begin{prop}
$g_{\lambda}(z)\rightarrow g_{\lambda}(x) $ a.e as $ y \rightarrow 0$ and
$$ g_{\lambda}(x)=\frac{1}{2}\lambda(h_{\lambda}\overline
{h_{\lambda}}g\chi_I)(x)+\frac{i}{2\pi}\lambda
h_{\lambda}(x)\ p.v\int_{I}\frac{(\overline {h_{\lambda}}g)(t)}{x-t}dt.$$

\end{prop}

\begin{dem}

We have : 
$$
g_{\lambda}(z)=\lambda h_{\lambda}(z)\frac{1}{2\pi i}\int_I \frac
{(\overline{h}_{\lambda}g)(t)}{t-z}dt \ \ \ \Im(z)>0.
$$

It is known that if $f$ is in $ L^{2}(\R) $, we can construct an
analytic function $ F $ on $\C_{+}$ by the integral formula :

$$
F(z)=\frac{1}{\pi i}\int_{\R}\frac{f(t)}{t-z}dt\ \ \ \mbox{and}\ \ \
F(z)\rightarrow f(x)+iHf(x)\ \qquad \mbox{a.e as}\ y\rightarrow 0
$$
where $ H $ is the Hilbert transform.\\
One can verify that the result remains valid if $f$ is a complex
valued function.\\
Let $f=\overline {h_{\lambda}}g\chi_{I}$, we see that :

$$
g_{\lambda}(z)=\lambda h_{\lambda}(z)\frac{1}{2\pi i}\int_I \frac
{f(t)}{t-z}dt \qquad \mbox{and as}\ y\rightarrow 0,
$$
$$
g_{\lambda}(z)\rightarrow \lambda
h_{\lambda}(x)\frac{1}{2}(f+iHf)(x) \qquad
\mbox{a.e}.
$$
Indeed, 
$$
g_{\lambda}(x)=\frac{1}{2}\lambda(h_{\lambda}\overline
{h_{\lambda}}g\chi_I)(x)+\frac{i}{2\pi}\lambda
h_{\lambda}(x)\ p.v\int_{I}\frac{(\overline {h_{\lambda}}g)(t)}{x-t}dt
\qquad \mbox{a.e}.
$$
\end{dem}
Usefull formulas of $g_{\lambda}$ for almost all x are :\\

$$
g_{\lambda}(x)=\left\{ \begin{array}{ll}
                     \frac{\lambda}{2(1+\lambda)}g(x)+\frac{i\lambda}{2\pi
                     (1+\lambda)}e^{iG_{1}(x)}\ p.v\int_I\frac{e^{-iG_{1}(t)}g(t)}{x-t}dt
                     & \mbox{for}\ x\in I\\       
                     \frac{\lambda}{2\pi
                     (1+\lambda)^{\frac{1}{2}}}e^{iG_{1}(x)}\ p.v 
\int_I\frac{e^{-iG_{1}(t)}g(t)}{x-t}dt  & \mbox{for}\ x\notin I    
\end{array} \right.
$$
\begin{rem}
\begin{enumerate}
\item Note that we do not need to integrate in the upper half plane as this
is the case in [7].
\item When $x\notin I$, observe that we can omit the p.v notation since
in this case $\displaystyle{p.v \int_I=\int_I}$.
\end{enumerate}
\end{rem}

\vspace{.5 cm}

\subsection{Nonsymmetrical case}

Let $I=[a,b]$ where $0 \leq a < b $, $g$ as in theorem 2.
\begin{prop}

$$h_{\lambda}(z)\rightarrow 
(1+\lambda\chi_{I}(x))^{\frac{-1}{2}} e^{iG_{2}(x)}\ \mbox{as} \ y\rightarrow 0 ,
\mbox{where} \ G_{2}(x)=\frac
{1}{2\pi}\ln(1+\lambda )\left[\frac{-1}{2}\ln\left
|\frac{1+b^2}{1+a^2}\right |+\ln\left | \frac {b-x}{a-x}\right|\right ].
$$

\end{prop}
\begin{dem}
What is changing here is the value of the integral
$\displaystyle{\int_{a}^{b}\frac{2t}{1+t^2}dt}$. Going back to the method of
proposition 1, we have 
$$ 
\int_{a}^{b}\frac{2t}{1+t^2}dt=\ln\left(\frac{1+b^2}{1+a^2}\right)
$$
and 
$$
h_{\lambda}(z)\rightarrow 
(1+\lambda\chi_{I}(x))^{\frac{-1}{2}}\displaystyle{ e^{\frac
{1}{2\pi}\ln(1+\lambda )\left[\frac{-1}{2}\ln\left
|\frac{1+b^2}{1+a^2}\right |+\ln\left | \frac {b-x}{a-x}\right|\right
]}}\ \mbox{as} \ y\rightarrow 0.
$$

Finally, 
$$
h_{\lambda}(x)=\left\{ \begin{array}{ll} 
                     (1+\lambda)^{\frac{-1}{2}}e^{iG_{2}(x)}
                     & \mbox{for}\ x\in I\\ 
                     e^{iG_{2}(x)}  & \mbox{for}\ x\notin I 
                     \end{array} \right.
$$
where
$$
 G_{2}(x)=\frac
{1}{2\pi}\ln(1+\lambda )\left[\frac{-1}{2}\ln\left
|\frac{1+b^2}{1+a^2}\right |+\ln\left | \frac {b-x}{a-x}\right|\right ].
$$
\end{dem}
For $g_{\lambda}$ we find the same expression, by substituting
$G_{1}(x)$ in place of $G_{2}(x)$.
\begin{prop}
$g_{\lambda}(z)\rightarrow g_{\lambda}(x) $ a.e as $ y \rightarrow 0$ and
for almost all x :
$$
g_{\lambda}(x)=\left\{ \begin{array}{ll}
                     \frac{\lambda}{2(1+\lambda)}g(x)+\frac{i\lambda}{2\pi
                     (1+\lambda)}e^{iG_{2}(x)}\ p.v\int_I\frac{e^{-iG_{2}(t)}g(t)}{x-t}dt
                     & \mbox{for} \ x\in I\\       
                     \frac{i\lambda}{2\pi
                     (1+\lambda)^{\frac{1}{2}}}e^{iG_{2}(x)}\  
\int_I\frac{e^{-iG_{2}(t)}g(t)}{x-t}dt  & \mbox{for}\ x\notin I    
\end{array} \right.
$$
\end{prop}

\vspace{.5 cm}

\section{Asymptotic behavior}

Suppose $I=]-a,a[$ and $x\notin\{-a,a\}$.\\ 
In this section we wish to estimate $g_\lambda(x)$ when $x\in
J=\R\setminus I$.\\
Let $\varphi (u)=a\displaystyle{\frac{1-e^{-u}}{1+e^{-u}}}$, a change of
variable and $\mathcal H$ be the following conditions :\\

 $$
\mathcal H \left\{ \begin{array}{ll}
               go\varphi\ \mbox{has no singularities in the strip }\
               0\leq \Im(z) < \pi \ \mbox{and does not accumulate on}\
               \Im(z)=\pi \ \ \mathcal C_{1}.\\
               \exists \delta\in\ ]0,1[, 
               \mbox{such that}\ go\varphi(z)
               =O(e^{\delta |\Re(z)|}) \\ \mbox{as}\ \Re(z)\rightarrow \pm
               \infty \ \mathcal \forall z\in \Gamma_j\ j\in\{1,2,3\}
               \ \ \qquad\qquad\qquad\qquad\ \ \mathcal C_{2}.  
\end{array} \right.
$$

$\Gamma=\displaystyle{\bigcup_j}\Gamma_j\cup [-R,R]$ is the following
contour and $\displaystyle b\in ]\pi,\frac{3\pi}{2}]$ :\\

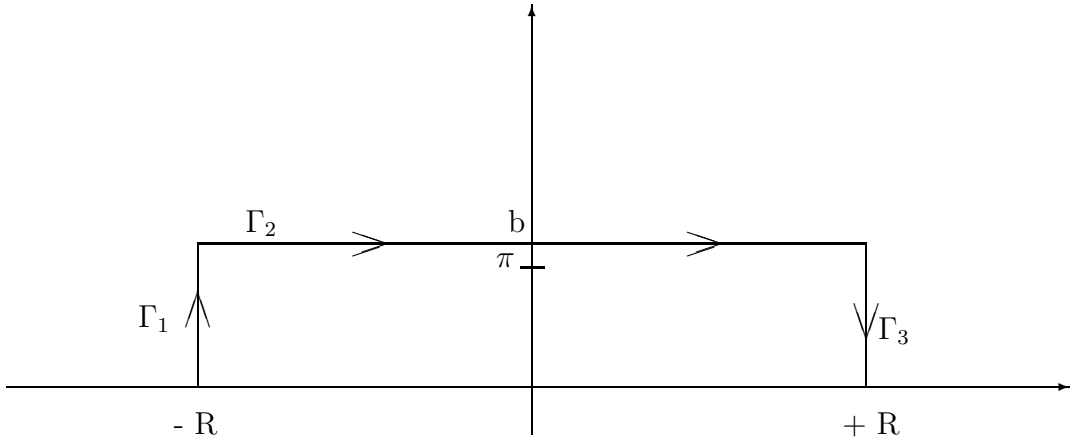
\input{contpress.latex}

\large {\bf Definition :}\\
$\mathcal {M}$ is the class of
functions verifying $\mathcal H$.

\begin{theoreme}
Let $I=]-a,a[, g\in \mathcal {M}$ we have :
$g_{\lambda}(x)=O(1)$.
\end{theoreme}
\begin{dem}
Let $g$ be in $\mathcal {M}$.\\
In order to simplify notations, suppose $go\varphi$ has no 
 singularities on $\Im (z)=\pi$.\\  
We have 
$$
g_{\lambda}(x)=\frac{\lambda}{2\pi (1+\lambda)^{\frac{1}{2}}}e^{iG_{1}(x)} 
\ \int_I\frac{e^{-iG_{1}(t)}g(t)}{x-t}dt,
$$
where 
$$
G_{1}(x)=\frac{1}{2\pi}\ln(1+\lambda )\ln\left | \frac {a-x}{a+x}\right |
$$

$$
\mbox{Put}\ \xi=\frac{\ln(1+\lambda)}{2\pi} ,\ u=-\ln\left | \frac
{a-t}{a+t}\right |\ \mbox{and} \ \alpha=\frac{a+x}{x-a}\  \mbox{so,}
$$
$$
\left |g_{\lambda}(x)\right|=\frac{\lambda a}{\pi
(1+\lambda)^{\frac{1}{2}}|a-x|}\left| \ \int_{\R}\frac{e^{i\xi u}e^{u}go\varphi(u)}
{(e^u+1)(e^u+\alpha)}du\right |.
$$
Note that :
$$
\frac{\lambda
a}{\pi(1+\lambda)^{\frac{1}{2}}(a-x)}=O(e^{\xi\pi}).\ \ \ \qquad (I)
$$ 
Denote $k(u,\xi)$, the quotient $\displaystyle{\frac{e^{i\xi u}e^{u}}{(e^u+1)(e^u+\alpha)}}$.\\
We use the residu theorem to estimate $\displaystyle\left |\int_{\R}k(u,\xi)
go\varphi(u)du\right |$.\\
Firstly, remark that singularities of the integral are of the form :
$$ 
i(\pi+2k\pi) \ \mbox{and}\ i(\pi+2k\pi)+\ln(\alpha)
$$
and
$$
\oint_{\Gamma} k(z,\xi)go\varphi(z)dz=\left (
Res(kg,i\pi)+Res(kg,i\pi+\ln(\alpha))\right )2\pi i.
$$
Calculations of both residus show that :

$$
Res(kg,i\pi)=\frac{e^{-\xi\pi}}{\alpha-1}go\varphi(i\pi)\ \ \mbox{and}
\ \ Res(kg,i\pi+\ln(\alpha))=\frac{e^{-\xi\pi}e^{i\xi\ln(\alpha)}}{1-\alpha}go\varphi(i\pi+\ln(\alpha)).
$$
On the other hand :
$$
\oint_{\Gamma} k(z,\xi)go\varphi(z)dz=\int_{-R}^{R}
k(u,\xi)go\varphi(u)du+\int_{0}^{b} k(R+iy,\xi)go\varphi(R+iy)dy$$
$$+\int_{R}^{-R} k(u+ib,\xi)go\varphi(u+ib)du+\int_{b}^{0} k(-R+iy,\xi)go\varphi(-R+iy)dy.
$$
When $R\rightarrow +\infty$ the first integral is the one we deal
with. If we set $I_{1},\ I_{2},\ I_{3}\ \mbox{and} \ I_{4}$ the four
integrals, we have :

$$
I_{1}=\frac{e^{-\xi\pi}}{\alpha-1}go\varphi(i\pi)+\frac{e^{-\xi\pi}e^{i\xi\ln(\alpha)}}{1-\alpha}go\varphi(i\pi+\ln(\alpha))-I_{2}-I_{3}-I_{4}.
$$
Put
$$
c_1(\alpha)=\frac{go\varphi(i\pi)}{\alpha-1} \ \mbox{and} \
c_2(\alpha,\xi)=\frac{e^{i\xi\ln(\alpha)}}{1-\alpha}go\varphi(i\pi+\ln(\alpha)),
$$ 
so
$$
\left|I_{1}\right|\leq
|c_1(\alpha)+c_2(\alpha,\xi)|e^{-\xi\pi}+|I_{2}|+|I_{3}|+|I_{4}|.
$$
Noting $c_2(\alpha,\xi)=c_2(\alpha)$, we get :
$$
\left|I_{1}\right|\leq
(|c_1(\alpha)|+|c_2(\alpha)|)e^{-\xi\pi}+|I_{2}|+|I_{3}|+|I_{4}|\ \ \ \
\ \ \ \ (II)
$$
\begin{rem}
Note that we can choose x such that $go\varphi$ does not vanish in
$\ln(\alpha)+i\pi$, so $c_2(\alpha,\xi) \neq 0$.
\end{rem} 
Moreover, $|I_2|=\displaystyle\int_{0}^{b} k(R+iy,\xi)go\varphi(R+iy)dy$ and $g$
satisfies $\mathcal C_{1}$, then \\$|go\varphi(R+iy)| \leq Me^{\delta R} \ \mbox{where} \ M >0$.
Thus,
$$
|I_2|\leq Me^{\delta R}\int_{0}^{b} |k(R+iy,\xi)|dy
$$
$$
\leq M\frac{e^{R}e^{\delta R}}{(e^R-1)(e^R-\alpha)}\int_{0}^{b}e^{-\xi y}dy
$$

$\mbox{when}\ R\rightarrow +\infty \ |I_{2}| \rightarrow 0.$\\
In a similar fashion we show that $\mbox{when}\ R \rightarrow +\infty \ |I_{4}| \rightarrow 0.$\\
For $I_{3}$, we see that 
$$
|I_{3}|\rightarrow \left |\int_{\R}k(u+ib,\xi)go\varphi(u+ib)du\right|\ 
\mbox{as} \ R\rightarrow+\infty.
$$

Without loss of generality, let $b=\displaystyle\frac{3\pi}{2}$ and if $\tilde{I}_{3}$\
denotes $\displaystyle\lim_{R\rightarrow +\infty}I_{3}$,
we get :
$$
|\tilde{I}_{3}|\leq
 e^{\frac{-3\pi}{2}\xi}\int_\R\frac{e^{\delta|t|}e^t}{|-ie^t+1||-ie^t+\alpha|}dt
$$
by the substitution $t=ln(r)$, we see that 
$$
|\tilde{I}_{3}|\leq e^{\frac{-3\pi}{2}\xi}\int_{0}^{+\infty}\frac{e^{|\ln(r^\delta)|}}{|-ir+1||-ir+\alpha|}dr
$$
$$\leq e^{\frac{-3\pi}{2}\xi}\left [
\int_{0}^{1}\frac{r^{-\delta}}{|-ir+1||-ir+\alpha|}dr+\int_{1}^{+\infty}\frac{r^{\delta}}{|-ir+1||-ir+\alpha|}dr
\right ]
$$
$$
\leq e^{\frac{-3\pi}{2}\xi}\left [
\int_{0}^{1}\frac{r^{-\delta}}{1+r^2}dr+\int_{1}^{+\infty}\frac{r^{\delta}}{1+r^2}dr\right ].
$$
Finally,
$$
|\tilde{I}_{3}|=O( e^{\frac{-3\pi}{2}\xi}).
$$
With (I) and (II) this ends the proof.
\end{dem}
\begin{rem}
\begin{enumerate}

\item Calculations show that the nonsymmetrical case give the same
results.
\item To bound $\tilde I_3$, the most important is $b\neq \pi$.
\item In the last proof, we supposed that $go\varphi$ has no
singularities on $\Im(z)=\pi$.\\
If not, $go\varphi$ has for instance $n$ pôles, $i\pi+\gamma_j$. Since
we can use the residu theorem, we get :\\
\begin{eqnarray*}
Res(kgo\varphi^,i\pi+\gamma_j) & = &
\frac{e^{i\xi(i\pi+\gamma_j)}e^{i\pi+\gamma_j}}{(e^{i\pi+\gamma_j}+1)(e^{i\pi+\gamma_j}+\alpha)}\lim_{u\rightarrow
i\pi+\gamma_j}go\varphi(u-i\pi-\gamma_j)^{m_j}\\
 & = & e^{-\xi\pi}c_j(\alpha,\xi,\gamma_j).
\end{eqnarray*}
$m_j$ is the ordre of the  p\^ole $i\pi+\gamma_j$. Finally :
$$
\oint_{\Gamma} k(z,\xi)go\varphi^(z)dz=2i\pi
e^{-\xi\pi}\sum_{j=1}^{n}c_j(\alpha,\xi,\gamma_j).
$$
Then we use the same methode to obtain the desired boundedness.

\end{enumerate}
\end{rem}

\begin{theoreme}
Suppose that $g$ satisfies $\mathcal C_{2}$  and $go\varphi$ is meromorphic in the strip 
$\Omega=\{0 < \Im(z) < \pi\}$  whose poles are in a finite
number in this open set.\\
Then  $g_{\lambda}(x)\rightarrow +\infty$ as $\lambda\rightarrow+\infty$.
\end{theoreme}

\begin{dem}
Employing the same methode as in the last theorem and without loss of
generality, assume that $go\varphi(z)$ has two poles, other than those of $k$, in $\Omega$, say
$\beta_{3}$ and $\beta_{4}$. We have,
$$
\oint_{\Gamma}k(z,\xi)g(z)dz=(Res(kg,i\pi)+Res(kg,i\pi+\ln(\alpha))+Res(kg,\beta_{3})+Res(kg,\beta_{4}))2\pi i
$$

but $Res(kg,i\pi)$ and $Res(kg,i\pi+\ln(\alpha))$ are known. If
$j\in\{3,4\}$ and $m_j$ is the order of the pole $\beta_j$, we see
that :

$$
Res(kg,\beta_j)=\frac{e^{i\xi\Re(\beta_j)}e^{-\xi\Im(\beta_j)}}{(1+e^{-\beta_j})
(1+\alpha e^{-\beta_j})}\lim_{z\rightarrow
\beta_j}go\varphi(z)(z-\beta_j)^{m_j},
$$
$$
Res(kg,\beta_j)=c_j(\beta_j,\xi)e^{-\xi\Im(\beta_j)}
$$
where,
$$
c_j(\beta_j,\xi)=\frac{e^{i\xi\Re(\beta_j)}}{(1+e^{-\beta_j})(1+\alpha e^{-\beta_j})}\lim_{z\rightarrow\beta_j}go\varphi(z)(z-\beta_j)^{m_j}
$$
If $|c_j(\beta_j,\xi)|=c_j(\beta_j)$, we see that
$|Res(kg,\beta_j)|=c_j(\beta_j)e^{-\xi\Im(\beta_j)}$.\\
On the other hand,
$$
\oint_{\Gamma}k(z,\xi)g(z)dz=I_1+I_2+I_3+I_4,
$$
therefore
$$
I_1=c_3(\beta_3,\xi)e^{-\xi\Im(\beta_3)}+c_4(\beta_4,\xi)e^{-\xi\Im(\beta_4)}
+[c_1(\alpha)+c_2(\alpha,\xi)]e^{-\xi\pi}-I_2-I_3-I_4. 
$$
Put $\displaystyle\tilde I_j=\lim_{R\rightarrow
+\infty}I_j$. Since $|\tilde I_2|=0$ and $|\tilde I_4|=0$ we have,
$$
|\tilde I_1|\geq e^{-\xi\pi}|c_1(\alpha)+c_2(\alpha,\xi)+c_3(\beta_3,\xi)e^{-\xi(\Im(\beta_3)-\pi)}+c_4(\beta_4,\xi)e^{-\xi(\Im(\beta_4)-\pi)}|-|\tilde I_3|.
$$
From the last theorem, $\exists M > 0 $ such that $-|\tilde I_{3}| \geq
-e^{-b\xi}M$, 
then
$$
e^{\xi\pi}|\tilde I_1|\geq |c_1(\alpha)+c_2(\alpha,\xi)+c_3(\beta_3,\xi)e^{-\xi(\Im(\beta_3)-\pi)}+c_4(\beta_4,\xi)e^{-\xi(\Im(\beta_4)-\pi)}|-e^{(\pi-b)\xi}M
$$
Note that $e^{(\pi-b)\xi}\rightarrow 0$ as $\xi\rightarrow +\infty$.
If for instance $\Im(\beta_3) \geq\Im(\beta_4)$ we deduce :
$$
e^{\xi\pi}|\tilde I_1|\geq e^{(\pi-\Im(\beta_3))\xi}|(c_1(\alpha)+c_2(\alpha,\xi))e^{(-\pi+\Im(\beta_3))\xi}+c_3(\beta_3,\xi)+c_4(\beta_4,\xi)e^{(\Im(\beta_4)-\Im(\beta_3))\xi}|.  
$$
Finally this last expression $\rightarrow +\infty$ as $\xi\rightarrow +\infty$.\\

\hspace{12cm} Q.E.D.
\end{dem}
\begin{coro}
If $go\varphi(z)$ verifies $\mathcal C_2$ then $\mathcal C_1$ is a necessary and
sufficient condition to obtain the boundedness of $g_\lambda(x)$.
\end{coro}

Since $x$ is arbitrarily taken, it can lie in a neighborhood $V$ of
positive measure, where the bounds of $I$ are not contained in $V$. We get the following :

\begin{coro}
If $go\varphi (z)$ verifies $\mathcal C_2$ and not $\mathcal C_1$ then
$g_\lambda (x) \rightarrow +\infty$ as $\lambda\rightarrow +\infty$
$\forall x \in V$. So we canot hope the ae convergence of $g_\lambda$
on $\T\setminus I$.
\end{coro}

Example 2 illustrates this.

\vspace{.5cm}
\section{Applications}

\vspace{.5cm}

Let $g(x)=g_1(x)+ig_2(x)$.
 
\begin{exmp}
$$
\left\{ \begin{array}{ll}
         g_1(x)=0 &\mbox{for}\ x\in]-\infty,-a[\\
         g_1(x)=\sqrt{a^2-x^2} &\mbox{for}\ x\in]-a,+a[\\
         g_1(x)=0 &\mbox{for} x\in]a,+\infty[
\end{array} \right.
$$
$$
\left\{ \begin{array}{ll}
         g_2(x)=-x-\sqrt{x^2-a^2} &\mbox{for}\ x\in]-\infty,-a[\\
         g_2(x)=-x &\mbox{for}\ x\in]-a,+a[\\
         g_2(x)=-x+\sqrt{x^2-a^2} &\mbox{for} x\in]a,+\infty[
\end{array} \right.
$$
In \ $]-a,+a[ , \ g(x)=\sqrt{a^2-x^2}-ix$.\\
To simplify calculations, assume $a=1$. We assert that $g \in \mathcal M$. In fact

$$
g_{1}o\varphi(t)=\sqrt{1-\left(\frac{1-e^{-t}}{1+e^{-t}}\right)^2}
=\frac{2e^{\frac{-t}{2}}}{1+e^{-t}}
$$
and
$$
g_{2}o\varphi(t)=\frac{e^{-t}-1}{e^{-t}+1}
$$
First, we are going to see  that $g_{1}o\varphi$ and  $g_{2}o\varphi$ satisfy
$\mathcal C_2$.\\

Observe that we have easily the equivalence :\\
  $go\varphi$ verifies $\mathcal C_2$ if and only if $g_{1}o\varphi$ and
$g_{2}o\varphi$  verify the same condition.\\

It is clear that boundedness  from above on $\Gamma_1$ and on
$\Gamma_2$ of $g_1$ and $g_2$ are sufficient to do the job.\\

\underline{Case $g_1$ on $\Gamma_1$} :\\

We have to find $M$ and $\delta\in ]0,1[$ such that :
$$
\left|\frac{2e^{-\left(\frac{-R+yi}{2}\right)}}{1+e^{R-yi}}\right|<Me^{\delta
R}
$$
where $R$ is some constant  $>0$ and $\forall y\in [0,b]$.\\
For this, one can for instance do 
$$
\left|\frac{2e^{-\left(\frac{-R+yi}{2}\right)}}{1+e^{R-yi}}\right|<\frac{2e^{\frac{R}{2}}}{1-e^R}<Me^{\delta R}.
$$

\underline{Case $g_1$ on $\Gamma_2$} :\\

As before we must find $M$ and $\delta\in ]0,1[$ such that :
$$
\left|\frac{2e^{-\left(\frac{ib+x}{2}\right)}}{1+e^{-(x+ib)}}\right|<Me^{\delta
|x|}\ \ \qquad \forall x\in\R.
$$
That is to say :
$$
\frac{2e^{-\frac{x}{2}}}{e^{-2x}+2\cos(b)e^{-x}+1}<e^{\delta|x|}.
$$
The first side, $h(x)$, is a positive function $C^\infty$
on $\R$ whose d\'enominator does not vanish. Furthermore,
$\displaystyle\lim_{x\rightarrow+\infty}h(x)=0$ and so does at
$-\infty$. Therefore there exists a maximum (function of $b$) such that
$h(x)<M(b)$,
where $\pi<b\leq\frac{3\pi}{2}$.\\

\underline{Case $g_2$ on $\Gamma_1$} :\\

It is easy to see that we can have $\delta\in ]0,1[$ and $M$ such that :

$$
\left|\frac{e^{-(-R+iy)}-1}{e^{-(-R+iy)+1}}\right|\leq Me^{-\delta R}.
$$

\underline{Case $g_2$ on $\Gamma_2$} :\\

A simple triangular inequatity application provides us the desired constants.

Note that the fact that  $go\varphi$ v\'erifies $\mathcal C_1$ is an
easy exercise.

\end{exmp}

\begin{exmp}
$$
g_1(x)=\frac{1}{1+x^2} \ \forall x\in \R
$$
$$
g_2(x)=\frac{-x}{1+x^2} \ \forall x\in \R
$$
After change of variable,
$$
go\varphi(t)=\frac{1}{2}\left (
\frac{(1+e^{-t})^2}{1+e^{-2t}}-i\frac{1-e^{-2t}}{1+e^{-2t}}\right )
$$
This example illustrates theorem 4. In fact,
$Res(kg,i\frac{\pi}{2})=e^{-\xi\frac{\pi}{2}}c(\alpha)$, 
where $c(\alpha)$ is a complex number.\\
We see that
$|g_\lambda(x)|\rightarrow +\infty \ \mbox{as} \ \lambda\rightarrow +\infty
$. 
\end{exmp}

\vspace{1cm}

\large{\Large\bf Aknowledjements} :\\

The subject was proposed to me by Professor P.Herve (I.U.T
d'énergétique et d'\'economie d'énergie. Ville d'Avray). I thank
J.Leblond (I.N.R.I.A-Sophia Antipolis) and Nicolas
Zakic (U.C.B.L)  for their valuable suggestions and
remarks.

}  
\end{document}

%% file: contpress.latex
\setlength{\unitlength}{0.012500in}%
\begingroup\makeatletter\ifx\SetFigFont\undefined
\def\x#1#2#3#4#5#6#7\relax{\def\x{#1#2#3#4#5#6}}%
\expandafter\x\fmtname xxxxxx\relax \def\y{splain}%
\ifx\x\y   
\gdef\SetFigFont#1#2#3{%
  \ifnum #1<17\tiny\else \ifnum #1<20\small\else
  \ifnum #1<24\normalsize\else \ifnum #1<29\large\else
  \ifnum #1<34\Large\else \ifnum #1<41\LARGE\else
     \huge\fi\fi\fi\fi\fi\fi
  \csname #3\endcsname}%
\else
\gdef\SetFigFont#1#2#3{\begingroup
  \count@#1\relax \ifnum 25<\count@\count@25\fi
  \def\x{\endgroup\@setsize\SetFigFont{#2pt}}%
  \expandafter\x
    \csname \romannumeral\the\count@ pt\expandafter\endcsname
    \csname @\romannumeral\the\count@ pt\endcsname
  \csname #3\endcsname}%
\fi
\fi\endgroup
\begin{picture}(445,180)(140,500)
\thinlines
\put(220,520){\line( 0, 1){ 60}}
\put(500,580){\line( 0,-1){ 60}}
\put(220,560){\line(-1,-3){  5}}
\put(220,560){\line( 1,-3){  5}}
\put(300,580){\line(-3, 1){ 15}}
\put(440,580){\line(-3, 1){ 15}}
\put(500,540){\line( 1, 3){  5}}
\put(300,580){\line(-3,-1){ 15}}
\put(440,580){\line(-3,-1){ 15}}
\put(140,520){\vector( 1, 0){445}}
\put(500,540){\line(-1, 3){  5}}
\put(240,585){\makebox(0,0)[lb]{\smash{\SetFigFont{12}{14.4}{rm}$\Gamma_2$}}}
\put(360,500){\vector( 0, 1){180}}
\put(220,580){\line( 1, 0){280}}
\put(355,570){\line( 1, 0){ 10}}
\put(490,500){\makebox(0,0)[lb]{\smash{\SetFigFont{12}{14.4}{rm}+ R}}}
\put(210,500){\makebox(0,0)[lb]{\smash{\SetFigFont{12}{14.4}{rm}- R}}}
\put(350,585){\makebox(0,0)[lb]{\smash{\SetFigFont{12}{14.4}{rm}b}}}
\put(345,570){\makebox(0,0)[lb]{\smash{\SetFigFont{12}{14.4}{rm}$\pi$}}}
\put(505,540){\makebox(0,0)[lb]{\smash{\SetFigFont{12}{14.4}{rm}$\Gamma_3$}}}
\put(195,545){\makebox(0,0)[lb]{\smash{\SetFigFont{12}{14.4}{rm}$\Gamma_1$}}}
\end{picture}